\documentclass[12pt]{amsart}

\usepackage{a4,color,graphics}

\usepackage{amscd}
\usepackage{amssymb}
\usepackage{amsmath}
\usepackage{hyperref}
\usepackage[all]{xy}
\xyoption{matrix}
\xyoption{arrow}

\newtheorem{theorem}{Theorem}[section]
\newtheorem{thm}[theorem]{Theorem}
\newtheorem{lem}[theorem]{Lemma}
\newtheorem{cor}[theorem]{Corollary}
\newtheorem{prop}[theorem]{Proposition}

\theoremstyle{definition}
\newtheorem{defn}[theorem]{Definition}
\newtheorem{eg}[theorem]{Example}

\newtheorem{quest}[theorem]{Question}

\theoremstyle{remark}

\numberwithin{equation}{section}

\def\vs#1{\vskip .#1 cm} 
\def\noi{\noindent}

\def\xrarrow{\xrightarrow} 
\def\xlarrow{\xleftarrow}

\def\smallcoprod{\,{\textstyle{\coprod}}\,}

\def\noteq{\neq}

\def\<{\left<}
\def\>{\right>}

\DeclareMathOperator{\Hom}{Hom}%

\newcommand{\field}[1]{\mathbb{#1}}
\newcommand{\ZZ}{\ensuremath{{\field{Z}}}}
\newcommand{\CC}{\ensuremath{{\field{C}}}}
\newcommand{\RR}{\ensuremath{{\field{R}}}}
\newcommand{\QQ}{\ensuremath{{\field{Q}}}}

\newcommand{\FF}{\ensuremath{{\field{F}}}}

\newcommand{\commentout}[1]{}

\def\el{\ell}

\newcommand{\cC}{\ensuremath{{\mathcal{C}}}}
\newcommand{\cD}{\ensuremath{{\mathcal{D}}}}

\newcommand{\cG}{\ensuremath{{\mathcal{G}}}}

\newcommand{\cI}{\ensuremath{{\mathcal{I}}}}

\newcommand{\cM}{\ensuremath{{\mathcal{M}}}}

\newcommand{\cP}{\ensuremath{{\mathcal{P}}}}

\newcommand{\cR}{\ensuremath{{\mathcal{R}}}}
\newcommand{\cS}{\ensuremath{{\mathcal{S}}}}

\def\a{\alpha}
\def\b{\beta}
\def\g{\gamma}

\def\d{\partial}

\def\del{\delete}

\def\f{\varphi}
\def\k{\kappa}
\def\r{\rho}
\def\s{\sigma}
\def\Sig{\Sigma}
\def\t{\tau}

\def\z{\zeta}

\def\delete{\backslash}

\def\ot{\leftarrow}
\def\then{\Rightarrow}

\def\onto{\twoheadrightarrow}
\def\ov{\overline}

\def\wt{\widetilde}
\def\what{\widehat}

\def\st{\,|\,}

\title{Families of Regular Matroids}
\author{Kiyoshi Igusa}
\address{Brandeis University}
\email{igusa@brandeis.edu}
\dedicatory{This paper is dedicated to my late friend, relative and fellow Mathematician\\ Dragan Acketa who loved matroids.}



\begin{document}

\maketitle

\setcounter{section}{0}
\setcounter{page}{1}


\begin{abstract}
This is an introductory paper about the category of regular oriented matroids (ROMs). We compare the homotopy types of the categories of regular and binary matroids. For example, in the unoriented case, they have the same fundamental group but we show that the higher homotopy groups are different for rank three regular and binary matroids. We also speculate on the possible impact of a recent theorem of Galatius \cite{Galatius:Free} computing the stable cohomology of the category of graphs and on possible applications to higher Reidemeister torsion.
\end{abstract}

\tableofcontents

%
%

\section*{Introduction}

This is an elementary account of the category of regular oriented matroids. ``Most'' regular matroids are either graph or cographic and, although we still don't know the homotopy type of these subcategories, they are the images of two categories of graphs which have the same stable cohomology, namely that of the infinite symmetric group. Thus it is plausible that the category of regular oriented matroids also has the same stable cohomology. I am also speculating that the category of homologically marked regular matroids has the same cohomology as the category of homologically marked graphs and therefore has higher torsion classes which will pull back to the Miller-Morita-Mumford or tautological classes for the Torelli group. (See \cite{Igusa:book} or \cite{I:ComplexTorsion} for more information about this.)

The paper is organized as follows. In the first section we review basic definitions of matroids and the homotopy types of categories. In section 2 we examine the category of binary matroids and regular matroids which form a full subcategory since every regular matroid has a unique binary structure. The first mod 2 homology functor induces a homotopy equivalence between the category $\cM(r,\FF_2)$ of binary matroids of rank $r$ and the category of $r$-dimensional vector spaces over the Galois field $\FF_2$. We conclude that $B\cM(r,\FF_2)\simeq BGL(r,\FF_2)$.

The category of regular matroids has a different homotopy type. We show that the homotopy fiber of the inclusion map $\cR(3,\FF_2)\to\cM(3,\FF_2)$ of the category of regular matroids of rank 3 to the category of binary matroids of the same rank is a wedge of eight 3-spheres. In section 3 we examine the category of regular oriented matroids. Since we allow change of sign as morphisms, we use an intrinsic sign-free definition of orientation. Again, the first homology functor gives a homotopy equivalence between the classifying space of the category of regular oriented matroids of rank $2$ with $BGL(2,\ZZ)$. Rank 3 will be discussed in a later paper.

In section 4 we discuss graphic and cographic matroids. The relation to graphs is mainly speculative. It is easy to see that the category of finite connected graphs and graph monomorphisms has the same classifying space as the symmetric group. By a recent result of S. Galatius \cite{Galatius:Free} the category of graphs and graph epimorphisms also has the same cohomology as the symmetric group in a stable range. These two categories map to the categories of cographic and graphic matroids respectively and both of these are full subcategories of the category of oriented regular matroids. We ask if all three of these matroid categories have the same cohomology in a stable range as the graph categories.

In section 5 we explain our motivation for studying regular matroids by reviewing the properties of higher Reidemeister torsion and the cohomology of the Torelli group and the category of homologically marked graphs. We then speculate that the higher Reidemeister torsion may exist as cohomology classes on the category of homologically marked regular oriented matroids. This would be a new and interesting method of analyzing these invariants.

Finally in section 6 we give the proof that the fiber poset $IR(r,\FF_2)$ for regular matroids of rank $r$ is simply connected for all $r$. In other words, it is the universal covering space for the category of regular matroids. Since $|IR(3,\FF_2)|$ is a wedge of 3-spheres, it is not 3-connected. This leaves open the question of whether $IR(r,\FF_2)$ could be 2-connected for all $r$.

I want to thank Keefe San Agustin for many helpful conversations and comments and corrections on many versions of this paper. I want to thank Milan Grulovic for his encouragement on this project. Finally I want to commemorate the life of Dragan Acketa. He was man of incredible energy and enthusiasm, whether it was rowing me and my wife across the Danube or talking to me about mathematics, extolling the virtues of matroids, I will always remember him and this paper is dedicated to him.


%
%

\section{Review of matroids and categories}

This is to go over the basic definitions and set the notation. We use standard notation for matroids and not quite standard notation for oriented matroids. The basic references are: \cite{Oxley}, \cite{WhiteBook87},\cite{BLvSWZ}.


\subsection{Matroids}

Matroids are finite sets together with the information of which subsets are \emph{independent}. This terminology suggests that the elements are vectors in a vector space and ``independence'' refers to ``linear independence.'' Minimal dependent subsets are called \emph{circuits}. This word suggests that we are thinking of the elements of the matroid as the edges of a graph and that ``independent'' sets are those that form a forest. Then the minimal dependent sets are cycles in the graph.

Usually, people talk about a ``matroid'' $M$ as a structure on a finite set which can be given in different ways.  To be more concrete we let $M$ denote the set of circuits of a matroid. Definition \ref{def:matroids by indep} and Proposition \ref{def:matroid by circuits} give two definitions of a matroid reflecting the vector analogy and the graph analogy.

\begin{defn}\label{def:matroids by indep}
A \emph{matroid} on a set $E$ is a collection of subsets called \emph{independent} sets satisfying the following properties.
\begin{enumerate}
\item The empty set is independent.
\item Any subset of an independent set is independent.
\item If $A,B\subseteq E$ are independent sets and $|A|>|B|$ then there is an element $a\in A\delete B$ so that $B\cup \{a\}$ is independent.
\end{enumerate}
We call this an \emph{independent subset matroid} on $E$.
Minimal dependent sets will be called \emph{circuits}. Maximal independent sets are called \emph{bases}. Complements of bases are called \emph{cobases}. Subsets of cobases are \emph{coindependent}. Subsets which are not coindependent are called \emph{covectors}. Minimal covectors are called \emph{cocircuits}.
\end{defn}

All bases have the same number of elements and this number is called the \emph{rank} of the matroid $M$. The coindependent sets satisfy the axioms for independent sets given above and therefore give a different matroid on the same set $E$. This is called the \emph{dual matroid} and denoted $M^\ast$. Its rank is $n-r$ where $n$ is the size of $E$ and $r$ is the rank of $M$.

\begin{prop}\label{def:matroid by circuits}
A collection of subsets $C\subseteq E$ forms the set of circuits of a matroid $M$ on $E$ if and only if they satisfy the following properties.
\begin{enumerate}
\item Circuits are nonempty.
\item If one circuit contains another then they are equal.
\item If $C_1,C_2$ are circuits and $e\in C_1\cap C_2$ then there is another circuit $C_3$ contained in $C_1\cup C_2$ so that $e\notin C_3$.
\end{enumerate}
\end{prop}

\begin{defn}
A \emph{representation} of a matroid on a set $E$ is a mapping
\[
    \r:E\to V
\]
of $E$ into a vector space $V$ over a field $K$ having the following properties.
\begin{enumerate}
\item A subset of $E$ is independent if and only if it maps monomorphically onto a linearly independent subset of $V$.
\item The image of $\r$ spans $V$.
\end{enumerate}
The matroid is then said to be \emph{representable} over the field $K$. Often we identify the elements of $E$ with their images in $V$. Thus we pretend that a representable matroid is a finite \emph{multiset} in $V$, i.e., a finite subset in which elements can have multiplicity.
\end{defn}

We need the following trivial fact about flats in a represented matroid. We recall that a subset $A$ of $E$ is called a \emph{flat} if any larger subset has strictly larger rank. Flats of corank 1 are called \emph{hyperplanes}.

\begin{prop}
Let $\r:E\to V$ be a representation of a matroid on $E$ over any field, and let $A$ be a subset of $E$.
\begin{enumerate}
\item The rank of $A$ is equal to the dimension of the linear span of $\r(A)\subseteq V$.
\item $A$ is a flat if and only if it is the inverse image of the linear span of its image in $V$.
\end{enumerate}
\end{prop}

There are two standard ways to give a matroid structure on a subset $E'$ of $E$. They are given by ``deleting'' and ``collapsing'' the elements of the complement of $E'$ in $E$. We will take the first definition.

\begin{defn} If $M$ is a matroid on $E$ and $E'$ is a subset of $E$ then the \emph{coinduced matroid} $M'=M|E'$ is given by either of the following equivalent conditions.
\begin{enumerate}
\item A subset of $E'$ is independent iff it is independent in $E$.
\item A subset of $E'$ is a circuit iff it is a circuit in $E$.
\end{enumerate}
We say that $M'$ is a \emph{submatroid} of $M$. The other standard matroid structure on $E'$ is the \emph{induced} matroid $(M^\ast|E')^\ast$.
\end{defn}

\begin{defn}
We define a \emph{morphism} of matroids
\[
    f:(E_1,M_1)\to (E_2,M_2)
\]
to be a monomorphism $f:E_1\to E_2$ so that a subset of $E_1$ is independent if and only if its image in $E_2$ is independent. Then $f$ gives an isomorphism
\[
    (E_1,M_1)\cong(f(E_1),M_2|f(E_1))
\]
Let $\cM$ denote the category of all matroids with these morphisms. Let $\cM(r)$ denote the full subcategory of matroids of rank $r$.
\end{defn}

We could ask the following.

\begin{quest}
What is the homotopy type of the category $\cM(r)$?
\end{quest}

However, we are more specifically interested in the homotopy type of the category of regular matroids. Since all regular matroids are binary, the category of regular matroids is a full subcategory of the category of binary matroids.



\subsection{Topology of categories}

We will need to assume that all categories are \emph{small}, which means the collection of objects is a set. This is not a problems since all of our categories are equivalent to small full subcategories and equivalent categories have the same homotopy type.

Next we review the construction of the classifying space $B\cC$ of a small category $\cC$. This is a topological space defined to be a quotient of a disjoint union of standard simplices
\[
    \Delta^k=\{(t_1,\cdots,t_k)\in\RR^k\st 1=t_0\ge t_1\ge t_2\ge \cdots\ge t_k\ge 0\}.
\]
For every $k\ge0$ we take one copy of $\Delta^k$ for every sequence of $k$ composable arrows in the category:
\[
    X_0\xlarrow{f_1}X_1\xlarrow{f_2}\cdots\xlarrow{f_k}X_k
\]
This gives a disjoint union
\[
    \coprod (f_1,\cdots,f_k)\times \Delta^k
\]
where the disjoint union is over all $k\ge0$ and all sequences of $k$ composable arrows $f_1,\cdots,f_k$. Each such sequence is considered to be one point (the first coordinate of the product) and is being used only to index multiple copies of the standard simplex $\Delta^k$.

The \emph{classifying space} $B\cC$ of $\cC$ is defined to be the quotient space
\[
    B\cC=\coprod (f_1,\cdots,f_k)\times \Delta^k/\sim
\]
where the equivalence relation is:
\[
    \left(
    X_{a(0)}\ot  \cdots\ot X_{a(k)}, t\in\Delta^k
    \right)\sim
        \left(
    X_0\ot \cdots\ot X_j, a_\ast(t)\in\Delta^j
    \right)
\]
for all nondecreasing mappings $a:[k]\to[j]$ where $[k]:=\{0,1,\cdots,k\}$. Here $a_\ast:\Delta^k\to\Delta^j$ is the simplicial map which is given by
\[
    a_\ast(t)_i=t_p
\]
where $p\in [k]$ is the smallest element so that $a(p)\ge i$.

Any functor $F:\cC\to \cD$ gives a continuous mapping
\[
    F_\ast:B\cC\to B\cD
\]
by
\[
    F_\ast(X_0\ot \cdots \ot X_k,t)= (FX_0\ot\cdots \ot FX_k,t)
\]

Some basic properties of this construction that we need are the following.

\begin{thm}
\begin{enumerate}\item If $\cC$ is a subcategory of $\cD$ then $B\cC$ is a closed subspace of $B\cD$.
\item A natural transformation of functors $\eta:F\to G$ gives a homotopy of induced maps on classifying spaces $F_\ast\simeq G_\ast$.
\item An equivalence of categories $F:\cC\cong \cD$ gives a homotopy equivalence of classifying spaces $F_\ast:B\cC\simeq B\cD$.
\item If $G$ is a discrete group and $\cC$ is the category with one object $\ast$ and $\Hom(\ast,\ast)=G$ then
\[
    B\cC =BG=K(G,1)
\]
is the classifying space of $G$ which is a path connected space with fundamental group $G$ and contractible universal covering space.
\item If $\cC$ is a connected groupoid (a category where all objects are isomorphic and all morphisms are isomorphism) then $B\cC\simeq BG$ where $G$ is the automorphism group of any object.
\end{enumerate}
\end{thm}

\begin{eg}\label{general linear groupoid} For any ring $R$ let $\cG\el(r,R)$ denote the category whose objects are right $R$-modules isomorphic to $R^r$ and whose morphisms are $R$-module isomorphisms. Then $\cG\el(r,R)$ is a connected groupoid. Therefore, as an example of the above theorem, we have a homotopy equivalence of classifying spaces
\[
    B\cG\el(r,R)\simeq BGL(r,R)
\]
where $GL(r,R)=Aut_R(R^r)$ is the group of invertible $r\times r$ matrices with coefficients in $R$.
\end{eg}

\begin{eg}\label{symmetric groupoid}
For any positive integer $n$, let $\cS(n)$ denote the category whose objects are finite sets with cardinality $n$ and whose morphisms are bijections of sets. Then $\cS(n)$ is a connected groupoid and the automorphism group of any object is the full permutation group of that set. Therefore,
\[
	B\cS(n)\simeq BS_n
\]
where $S_n$ is the symmetric group on $n$ letters.
\end{eg}


%
%

\section{Binary matroids}

We now examine the space of binary matroids of a fixed rank. We recall that a matroid is called \emph{binary} if and only if it is representable over $\FF_2$, the field with 2 elements. For example, take the nonzero element of $\FF_2$ $n$ times. This is a binary matroid with $n$ elements in which the entire set forms the unique cocircuit and any two elements form a circuit. This is the ``uniform matroid'' $U(1,n)$.

We will use one well-known theorem which characterizes binary matroids.

\begin{thm}\label{thm:binary matroid even intersection characterization}
A matroid is binary if and only if the intersection of any circuit and cocircuit has an even number of elements.
\end{thm}

We also need the \emph{rigidity} of binary matroids:

\begin{thm}\label{thm:rigidity of binary matroids}
Suppose that $\r:E\to V,\r':E\to V'$ are two $\FF_2$-representations of a binary matroid on $E$. Then there is a unique linear isomorphism $\f:V\to V'$ so that $\r'=\f\circ\r$.
\end{thm}

\begin{proof} Suppose that $\r:E\to V$ is a representation over any field $K$.
Choose a basis $B=\{b_1,\cdots,b_r\}$ for $E$. Then $\r(B)$ is a basis for $V$. This gives an isomorphism $V\cong K^r$. For any element $e\in E\delete B$, there is a unique circuit $C_e$ in $E$ so that $e\in C_e\subseteq B\cup\{e\}$. This is called a \emph{basic circuit}. The $i$th coordinate of $\r(e)$ will be nonzero if and only if $b_i$ is in the basic circuit $C_e$. When $K=\FF_2$, this determines $\r(e)$ uniquely for every $e\in E\delete B$.
\end{proof}

We will see that this theorem implies that the space of binary matroids of rank $r$ is the classifying space for the group $GL(r,\FF_2)$. More precisely, the category $\cM(r,\FF_2)$ of binary matroids of rank $r$ is homotopy equivalent to the category $\cG \el(r,\FF_2)$ of $r$-dimensional vector spaces over the Galois field $\FF_2$ and this homotopy equivalence
\[
    \cM(r,\FF_2)\xrarrow{\simeq} \cG \el(r,\FF_2)
\]
is given by the mod-2 homology functor as we will explain below.



\subsection{Homology modulo 2}

If $M$ is a matroid on $E$ we define the \emph{first homology group} $H_1(M)=H_1(E, M)$ to be the additive group generated by the elements of $E$ modulo the relation
\[
\sum_{e\in C}e=0
\]
for all circuits $C$ in $M$. For any abelian group $A$ we define homology with coefficients in $A$ as the tensor product: $H_1(M;A):=H_1(M)\otimes A$.

It is easy to see that the image of any basis in $E$ will generate $H_1(M)$. Therefore, the rank of $H_1(M)$ over any field is bounded by the rank of $M$.

We define the \emph{first cohomology} of $M$ with coefficient in $A$ to be $H^1(M;A)=\Hom(H_1(M),A)$. If $A=F$ is a field then $H^1(M;F)$ is just the vector space of all mappings
\[
    f:E\to F
\]
which sends every circuit to zero in the sense that $f(C)=\sum_{e\in C} f(e)=0$. We call such functions \emph{cocycles} with coefficients in $F$. In the special case $F=\FF_2$, we often identify a cocycle $f$ with its \emph{support} $\{e\in E\st f(e)\neq0\}=f^{-1}(1)$.

It is easy to see that $H^1(M;F)\cong \Hom_F(H_1(M;F),F)$. So, the rank of $H_1(M;F)$ is equal to the rank of $H^1(M;F)$ for any field $F$. Also, we have a \emph{universal cocycle}
\[
    \iota: E\to H_1(M;\FF_2)
\]
which is the identity map in $\Hom_{\FF_2}(H_1(M;\FF_2),H_1(M;\FF_2))\cong H^1(M,H_1(M;\FF_2))$

Theorem \ref{thm:binary matroid even intersection characterization} has the following interpretation in terms of mod-2 cohomology.

\begin{thm} If $M$ is a matroid of rank $r$ on a set with $n$ elements then the following are equivalent.
\begin{enumerate}
\item $M$ is binary.
\item $H^1(M;\FF_2)\cong \FF_2^r$.
\item Every cocircuit of $M$ is (the support of) a cocycle modulo 2.
\item The cocircuits of $M$ are the mod-2 cocycles with minimal support.
\item $H^1(M^\ast;\FF_2)\cong \FF_2^{n-r}$ (i.e., $M^\ast$ is binary).
\end{enumerate}
\end{thm}

By the rigidity of binary matroids (Theorem \ref{thm:rigidity of binary matroids}), this theorem implies the following.

\begin{cor} The universal cocycle $
    \iota:E\to H_1(M;\FF_2)
$
for a binary matroid $M$ on $E$ is a \emph{universal representation} in the sense that any $\FF_2$-representation of $M$ factors uniquely through $\iota$.
\end{cor}



\subsection{Category of binary matroids}

Let $\cM(r,\FF_2)$ be the category of binary matroids of rank $r$. We will calculate the homotopy type of the classifying space $B\cM(r,\FF_2)$. First, we need the following definition.

Suppose that $F$ is any field and
\[
    \r_i:E_i\to V
\]
are representations of two matroids $M_1,M_2$ on $E_1,E_2$ with the same target space $V\cong F^r$. Then we define the union
\[
    M_1\cup_V M_2
\]
to be the representable matroid structure on $E_1\coprod E_2$ defined by the representation
\[
    \r_1\smallcoprod \r_2:E_1\smallcoprod E_2\to V.
\]
The key point is that $M_1,M_2$ are submatroids of $M_1\cup_VM_2$ of full rank. When $F=\FF_2$, $M_1,M_2$ and $M_1\cup_VM_2$ are all binary matroids.

\begin{thm}
The mod-2 homology functor 
\[
	F:= H_1(-;\FF_2):\cM(r,\FF_2)\to \cG\el(r,\FF_2)
\]
induces a homotopy equivalence on classifying spaces
\[
    F_\ast:B\cM(r,\FF_2)\simeq B\cG\el(r,\FF_2)\simeq BGL(r,\FF_2).
\]
\end{thm}

\begin{proof}
A (right) inverse functor
\[
    G:\cG\el(r,\FF_2)\to \cM(r,\FF_2)
\]
is given by sending the $\FF_2$-vector space $V$ to the underlying set $V$ with the given linear matroid structure. By the rigidity of binary matroids (Theorem \ref{thm:rigidity of binary matroids}), this is a section of the functor $F$. In other words,
\[
    H_1(G(V);\FF_2)\cong V
\]
giving a homotopy
\[
    F\circ G\simeq id_{\cG\el(r,\FF_2)}.
\]
Conversely, we have a homotopy
\[
    G\circ F\simeq id_{\cM(r,\FF_2)}
\]
given by the sequence of natural transformations
\[
    GF(M)=H_1(M;\FF_2)\subseteq GF(M)\cup_V M\supseteq M
\]
where $V=H_1(M;\FF_2)$ and the union $GF(M)\cup_V M$ is over the universal representation
\[
    \iota:E\to H_1(M;\FF_2)=V
\]
and the identity mapping $H_1(M;\FF_2)\to V$.
\end{proof}



\subsection{Category of regular matroids}

In lieu of the definition we use the following well-known characterization of regular matroids due to Tutte with terms defined below.

\begin{thm}[Tutte] Every regular matroid is binary and
a binary matroid is regular if and only if it has no minor isomorphic to $F_7$ or $F_7^\ast$.
\end{thm}


Let $\cR(r,\FF_2)$ denote the category of regular matroids of rank $r$. By Tutte's theorem above, this is the full subcategory of the category $\cM(r,\FF_2)$ whose objects are those binary matroids of rank $r$ which have no minor isomorphic to the Fano matroid $F_7$ or its dual $F_7^\ast$.

\begin{defn}
The \emph{Fano matroid} $F_7$ is defined to be the binary matroid given by the seven nonzero elements of $\FF_2^3$:
\[
    F_7=\left[
    \begin{array}{ccc|cccc}
    1 & 0 & 0 & 1 & 1 & 0 & 1\\
    0 & 1 & 0 & 1 & 0 & 1 & 1\\
    0 & 0 & 1 & 0 & 1 & 1 & 1
    \end{array}
    \right]
\]
\end{defn}

If this is abbreviated $F_7=[I_3|C]$ then the dual is
\[
    F_7^\ast=[-C^t|I_4]=\left[
    \begin{array}{ccc|cccc}
    1 & 1 & 0 & 1 & 0 & 0 & 0\\
    1 & 0 & 1 & 0 & 1 & 0 & 0\\
    0 & 1 & 1 & 0 & 0 & 1 & 0\\
    1 & 1 & 1 & 0 & 0 & 0 & 1
    \end{array}
    \right]    
\]
It is very easy to see that this standard construction gives the dual of a binary matroid. The calculation
\[
    [-C^t|I_4]\left[
    \begin{matrix}
    I_3\\
    C^t
    \end{matrix}
    \right]=C^t-C^t=0
\]
shows that the dual of the row space of each of these matrices is equal to the null space of the other matrix. But the elements with minimal support in the null space of a matrix over $\FF_2$ are the cocircuits of the binary matroid. So, the cocircuits of each matroid are the circuits of the other matroid. So, they are dual matroids.

To say that a binary matroid on $E$ has a \emph{minor} isomorphic to $F_7$ means that it has seven elements $e_1,\cdots,e_7$ and three cocycles $c_1,c_2,c_3$ so that
\[
    [c_i(e_j)]=[I_3|C]
\]
and so that the kernel $\{x\in E\st c_1(x)=c_2(x)=c_3(x)=0\}$ is a codimension 3 flat in $E$. Similarly, $F_7^\ast$ is a minor if there are seven elements in $E$ and a codimension 4 flat in $E$ which is the kernel of four cocycles satisfying the analogous equation.

Some immediate consequences of Tutte's theorem are the following.

\begin{prop}\begin{enumerate}
\item Any binary matroid of rank $\le 2$ is regular.
\item Any binary matroid of rank $r$ on a set with $\le r+2$ elements is regular.
\end{enumerate}
\end{prop}

\begin{cor}\label{cor:R=M for r le 2}
$\cR(r,\FF_2)=\cM(r,\FF_2)\simeq \cG\el(r,\FF_2)$ for $r\le2$.
\end{cor}

For $r\ge3$ this is not true because $F_7\in\cM(3,\FF_2)$.

Since $\cM(r,\FF_2)\simeq\cG\el(r,\FF_2)$, the difference between $\cR(r,\FF_2)$ and $\cM(r,\FF_2)$ is given by the homotopy fiber of the mod 2 homology functor
\[
    H_1(-,\FF_2):\cR(r,\FF_2)\to \cG\el(r,\FF_2).
\]

We recall that for any functor $F:\cC\to\cG$ from a small category $\cC$ to a connected groupoid $\cG$ we can construct the \emph{fiber category} $\wt\cC$ as follows. Choose one object $Z$ of $\cG$. Then the objects of $\wt\cC$ will be pairs $(X,f)$ where $X$ is an object of $\cC$ and $f:Z\to FX$ is an isomorphism. A morphism $(X,f)\to (Y,g)$ is defined to be a morphism $h:X\to Y$ so that $g=Fh\circ f:Z\to FX\to FY$. The group $G=Aut_\cG(Z)$ acts on $\wt\cC$ on the right by $(X,f)g=(X,fg)$. This is a free action since $(X,f)$ is never equal to $(X,fg)$ unless $g=1=id_Z$. Therefore, $G$ acts freely on the topological space $B\wt\cC$ and the quotient is homeomorphic to $B\cC$. Since these spaces are CW-complexes and the action of $G$ on $B\wt\cC$ is cellular, it is a locally trivial regular covering space which is classified by the mapping
\[
    F_\ast:B\cC\to B\cG\simeq BG.
\]
This is an example of Quillen's Theorem B \cite{Quillen73}.

\begin{defn}\label{def:fiber poset for regular binary matroids}
For $r\ge 1$ let $IR(r,\FF_2)$ denote the poset of all subsets of $\FF_2^r\backslash 0$ which are regular matroids of rank $r$ ordered by inclusion.
\end{defn}

We recall that the \emph{geometric realization} $|P|$ of a poset $P$ is the geometric realization of the simplicial complex whose $k$ simplices are the $k+1$ chains in $P$ which are the sequences $x_0<x_1<\cdots<x_k$ in $P$. This is also homeomorphic to the classifying space of $P$ considered as a category. If $P$ is a simplicial complex then $|P|$ is (the geometric realization of) the first barycentric subdivision of $P$.

\begin{prop}\label{prop: fiber poset for regular matroids}
$|IR(r,\FF_2)|$ is homotopy equivalent to the classifying space of the fiber category of $H_1(-,\FF_2):\cR(r,\FF_2)\to \cG\el(r,\FF_2)$. In other words, we have a homotopy fiber sequence:
\[
    |IR(r,\FF_2)|\to B\cR(r,\FF_2)\to B\cG\el(r,\FF_2)
\]
\end{prop}

\begin{proof}
The fiber category of $\cR(r,\FF_2)$ over $\FF_2^r\in \cG\el(r,\FF_2)$ consists of all pairs $((E,M),\r)$ where $M$ is a regular matroid on $E$ of rank $r$ and $\r:E\to\FF_2^r$ is a representation of $M$. (By rigidity of binary matroids, a representation over $\FF_2$ of any regular matroid is the same as an isomorphism $H_1(M,\FF_2)\cong \FF_2^r$.) The poset $IR(r,\FF_2)$ is the full subcategory of all such pairs where $\r:E\to\FF_2^r$ is an inclusion map. And a deformation retraction of the fiber category $\widetilde\cR(r,\FF_2)$ to $IR(r,\FF_2)$ is given by mapping $((E,M),\r)$ to the image $\r(E)$ minus 0.
\end{proof}

From the above proof we extract the following definitions.

\begin{defn}\label{mod2 homologically marking}
An isomorphism $H_1(M;\FF_2)\cong \FF_2^r$ will be called a \emph{mod-2 homological marking} for the regular matroid $M$. The fiber category $\wt\cR(r,\FF_2)$ is thus the category of mod-2 homologically marked regular matroids of rank $r$. We will refer to this category later as $\cI\cR(r,\FF_2)$.
\end{defn}



\subsection{Homotopy type of the fiber poset}

We will go over some basic homotopy properties of the poset $IR(r,\FF_2)$.

\begin{thm}\label{thm:IR(F2) is simply connected}
$IR(r,\FF_2)$ is simply connected for all $r$.
\end{thm}

When we say that a poset or category is $k$-connected we mean that its geometric realization is $k$-connected. Also the homotopy groups of a poset or category is defined to be those of its geometric realization. We will give the proof of this theorem in the last section \ref{last section}. By the fiber sequence in Proposition \ref{prop: fiber poset for regular matroids} above, this gives the following.

\begin{cor}
For all $r$, the functor $H_1(-;\FF_2)$ induces an isomorphism
\[
    \pi_1 \cR(r,\FF_2)\cong GL(r,\FF_2).
\]
Furthermore, the universal covering space of $B\cR(r,\FF_2)$ is homotopy equivalent to $|IR(r,\FF_2)|$.
\end{cor}

For $r=1,2$ we know by Corollary \ref{cor:R=M for r le 2} that $IR(r,\FF_2)$ is contractible. However, this is not true for $r=3$. To compute the homotopy type of $IR(3,\FF_2)$ we use shellable simplicial complexes which are defined as follows. (Shellable complexes are defined more generally. See, e.g., \cite{BLvSWZ}. See \cite{WhiteBook92}, Chap 7 for more details about shellable simplicial complexes.)

\begin{defn}
The empty simplicial complex is defined to be a shellable complex of any dimension and for any $n\ge0$, a finite $n$-dimensional simplicial complex $K$ is \emph{shellable} if
\begin{enumerate}
\item $K$ is the union of closed $n$-simplices $S_1,\cdots,S_k$ and
\item $(S_1\cup\cdots\cup S_j)\cap S_{j+1}$ is a shellable $n-1$ dimensional complex for all $j$.
\end{enumerate}
The sequence of $n$-simplices $S_1,\cdots,S_k$ is called a \emph{shelling} of $K$.
\end{defn}

Using the obvious fact that any union of $n-1$ dimensional faces of $\Delta^n$ is either empty, contractible, or all of $\d\Delta^n$, we see that any connected shellable $n$-dimensional simplicial complex is homotopy equivalent to a wedge of $n$-spheres.

\begin{thm}
$IR(3,\FF_2)$ is homotopy equivalent to a wedge of eight $3$-spheres:
\[
    |IR(3,\FF_2)|\simeq\bigvee_8 S^3.
\]
\end{thm}

\begin{proof}
We will show that $|IR(3,\FF_2)|$ is the first barycentric subdivision of a $3$-dimensional simplicial complex which is connected and shellable with Euler characteristic -7. The theorem follows.

Recall that $F_7$ is the matroid consisting of all 7 of the nonzero elements of $\FF_2^3$. The poset $IR(3,\FF_2)$ is the poset of all proper subsets of this set which span $\FF_2^3$. There is an order reversing bijection of this poset with the poset $P$ of all complements of these sets. These are the nonempty subsets of the seven point set $S=\FF_2^3\backslash 0$ which have at most 4 elements and are not equal to one of the seven 4-circuits which are the complements of the seven hyperplanes. For any $\s\in P$, $\s$ is a subset of $S$ and every nonempty subset of $\s$ is also an element of $P$. Therefore $P$ is a simplicial subcomplex of the 6-simplex $\Delta^S$.

The Euler characteristic of $P$ is easy to compute:
\[
    \chi(P)=\binom71-\binom72+\binom73-\binom74+7=-7
\]
The shellability of $P$ follows from the following lemma and the easy fact that the intersection of any two missing tetrahedra is an edge.
\end{proof}

\begin{lem}
Suppose that $K$ is a $d$-dimensional subcomplex of the standard $n$-simplex $\Delta^n$ so that $K$ contains the $d-1$ skeleton of $\Delta^n$ and any two missing $d$-simplices intersect in a $d-2$ simplex. Then $K$ is shellable.
\end{lem}

I don't know who noticed this fact first.

\begin{proof}
If $d=0$ then $K$ is shellable by definition. If $n=d$ then $K=\Delta^d$ is shellable. So suppose $0<d<n$ and proceed by double induction on the pair $(d,n)$.

Since $K$ has dimension $d$, it contains at least one $d$-simplex. Since $n>d$ there is a vertex $v$ in $\Delta^n$ which is not in this $d$-simplex. It is easy to see that the link $L(v)$ of $v$ in $K$ is a $d-1$ dimensional complex satisfying the conditions of the lemma and is therefore shellable by induction on $d$. The intersection $K_0$ of $K$ with the $n-1$ face $\d_v\Delta^n$ of $\Delta^n$ opposite the vertex $v$ is also shellable by induction on $n$. But $K_0$ contains the entire $d-1$ skeleton of $\d_v\Delta^n$. So, if we take any shelling of $K_0$ and follow it with the cone of any shelling of $L(v)$ we get a shelling of $K=K_0\cup_{L(v)}CL(v)$.

To verify this last statement note that each $d-1$ simplex $\s$ in a shelling of $L(v)$ is attached on a shellable $d-2$ subcomplex $M\subseteq \d\s$. In the corresponding shelling of $K$, the cone of $\s$ is attached along the union of the cone $CM$ of $M$ with $\s$. The required shelling of $CM\cup_M\s$ can be given by first taking the cone of a shelling of $M$ and then add $\s$ last along $M$.
\end{proof}

The poset $IR(3,\FF_2)$ has an action of $G=SL(3,\FF_2)$ which makes the homology of $IR(3,\FF_2)$ into a representation of $G$. The order of this well-known finite simple group is $|G|=2^3\cdot3\cdot7=168$. So you can write down all the elements and see that they have order 1,2,3,4,7 and there are two conjugacy classes of elements of order 7. To see what they are, use $\FF_8$ as a model for $\FF_2^3$. The group $\FF_8^\times$ has 6 elements of order 7. They are $\a,\a^2,\a^4$ which are the roots of the polynomial $x^3+x+1$ and $\a^3=\b,\b^2,\b^4$ which are the three roots of $x^3+x^2+1$. Multiplication by $\a,\b$ give elements of $SL(3,\FF_2)$ which are not conjugate since they have different characteristic polynomial. With this observation it is easy to compute the character table of $SL(3,\FF_2)$:

\begin{table}[htdp]
\caption{$SL(3,\FF_2)$, $z=\g+\g^2+\g^4$, $\g=e^{2\pi i/7}$}
\begin{center}
\begin{tabular}{c|cccccc}
$C_i$    & 1 & 2 & 3 & 4 & 7($\a$) & 7($\b$)\\
$|C_i|$  & 1 & 21 & 56 & 42 & 24 & 24\\
    \hline
    $\chi_0$ & 1 & 1 & 1 & 1 & 1 & 1 \\
    $\chi_1$ & 6 & 2 & 0 & 0 & -1 & -1\\
    $\chi_2$& 7 & -1 & 1 & -1 & 0 & 0\\
    $\chi_3$ & 8 & 0 & -1 & 0 & 1 & 1\\
    $\chi_4$ & 3 & -1 & 0 & 1 & $z$ & $\ov z$\\
    $\chi_5$ & 3 & -1 & 0 & 1 & $\ov z$ & $z$    
\end{tabular}
\end{center}
\label{character table of SL(3,F2)}
\end{table}%

\begin{thm}
$H_3(IR(3,\FF_2),\CC)$ is the unique $8$ dimensional irreducible representation of $SL(3,\FF_2)$.
\end{thm}

\begin{proof}
As we saw above, $|IR(3,\FF_2)|$ is the first barycentric subdivision of the $3$-dimensional simplicial complex given by nonempty subsets of $\FF_3^\times$ with $\le 4 $ elements which are not $4$-circuits. The action of $G=SL(3,\FF_2)$ on the set of oriented $4$-element subsets of $\FF_8^\times$ has two orbits, one is the set of $4$-circuits. (Here the \emph{orientation} of a set is an ordering up to even permutation.) Thus, as an element of $K_0\CC G$, we have:
\[
    \sum (-1)^i [H_i(IR(3,\FF_2),\CC)]=[C_0]-[C_1]+[C_2]-[C_3^{(b)}]
\] where $C_k$ is the vector space over $\CC$ of generated by oriented $k+1$ element subsets of $\FF_8^\times$ considered at the $k+1$ fold exterior power of $C_0$. $C_3=C_3^{(a)}\oplus C_3^{(b)}$ where $C_3^{(a)}$ is the vector space spanned by oriented $4$-circuits in $\FF_8^\times$.

The group $G=SL(3,\FF_2)$ acts by even permutations on the 7 elements of $\FF_8^\times$ since otherwise it would have an index 2 subgroup which is impossible. Therefore, the taking of complements is an isomorphism of $G$-modules $C_k\cong C_{5-k}$ and $C_3^{(a)}\cong C_2^{(a)}$, the vector space of oriented hyperplanes in $\FF_2^3$. So,
\[
    \sum (-1)^i [H_i(IR(3,\FF_2),\CC)]=[C_0]-[C_1]+[C_2^{(a)}]
\]
Computing the characters of these three representations is easy:
\[
    (7,3,1,1,0,0)-(21,1,0,-1,0,0)+(7,-1,1,-1,0,0)=(-7,1,2,1,0,0).
\]
Subtracting the trivial character (for $H_0(IR(3,\FF_2))$) and changing sign, we get the theorem.
\end{proof}

\begin{cor}
The classifying space for the category of regular matroids of rank $3$ has the rational homotopy type of a point.
\end{cor}

\begin{proof}
The Serre spectral sequence for the homotopy fiber sequence (Proposition \ref{prop: fiber poset for regular matroids}) collapses making $B\cR(3,\FF_2)$ rationally homotopy equivalent to $B\cG\el(3,\FF_2)$. This in turn is rationally acyclic since $GL(3,\FF_2)$ is a finite group.
\end{proof}


%
%

\section{Freely oriented matroids}

Now we move onto oriented matroids. We want every \emph{unoriented} graph with $n$ edges to give an oriented matroid on a set with $n$ elements. More precisely we get $2^n$ equivalent oriented matroids on this set. We want to form a category whose objects are equivalence classes of oriented matroids. To make this more concrete we will replace $E\times\{+,-\}$ with a \emph{two-fold covering} $\wt E$ of $E$. This means we are given an epimorphism
\[
    p:\wt E\to E
\]
so that the inverse image of every element of $E$ is two elements of $\wt E$. The free involution on $\wt E$ which switches all such pairs will be denoted with an overline: $e\mapsto\ov e$.


\subsection{Basic definitions}

\begin{defn} A \emph{freely oriented matroid} on $E$ is defined to be a pair $(\wt E,\wt M)$ where $\wt E$ is a two-fold covering of $E$ and $\wt M$ is a collection of nonempty subsets $C$ of $\widetilde{E}$ called \emph{signed circuits} satisfying the following conditions.
\begin{enumerate}
\item $C$ maps monomorphically into $E$, i.e., $C\cap\ov{C}=\emptyset$.
\item If $C$ is a signed circuit then so is $\ov{C}$.
\item If $C_1\subseteq C_2\cup\ov{C_2}$ then either $C_1=C_2$ or $C_1=\ov{C_2}$.
\item Suppose that $C_1\neq{C_2}$ and $e\in C_1\cap{C_2}$. Then there is a third signed circuit $C_3$ contained in $C_1\delete e\cup \ov{C_2\delete{e}}$.
\end{enumerate}
\end{defn}

\begin{eg}
An unoriented graph $G=(G_0,G_1)$ gives a freely oriented matroid on its set of half-edges $\wt G_1$. If we view the two half-edges which make up a single edge as being the two possible orientations of that edge, then a circuit $C\subseteq \wt G_1$ is a minimal oriented cycle. Such matroids are called \emph{graphic matroids}.
\end{eg}

\begin{eg}
Suppose that $\Delta_n$ is a Dynkin diagram with $n$ vertices. Then the set of positive roots $\Phi_+$ is a finite spanning set of vectors in $\RR^n$ with a two fold covering given by the set of all roots $\Phi(\Delta_n)=\Phi_+\cup\Phi_-$. The inclusion map $\Phi(\Delta_n)\to \RR^n$ is a representation of an oriented matroid. It is not too hard to see that this matroid is graphic if and only if the Dynkin diagram is of type $A_n$. (Otherwise, there would be two independent elements $\a,\b\in\Phi$ so that $\a+\b,\a+2\b\in\Phi$, which is impossible in any regular matroid.)
\end{eg}

It is clear that we can transform any freely oriented matroid into an oriented matroid by choosing a $\ZZ/2$ equivariant mapping
\[
    \widetilde{E}\to\{+,-\}.
\]
This is equivalent to choosing a section of the covering map $\wt E\to E$. Either of these equivalent choices will be called an \emph{orientation} of $\wt E$.
This implies that all results about oriented matroids apply, with suitable change in terminology, to freely oriented matroids. For example:

\begin{prop}
The images in $E$ of the signed circuits in $\wt E$ form a circuit matroid $M=p(\wt M)$ on $E$.
\end{prop}

The \emph{signed cocircuits} are the subsets $C^\ast$ of $\wt E$ which
\begin{enumerate}
\item map monomorphically onto a cocircuit in $E$ and
\item meet any signed circuit $C$ if and only if they meet $\ov C$. I.e, for any circuit $C\subseteq\wt E$, $C\cap C^\ast\noteq\emptyset$ if and only if $\ov C\cap C^\ast\noteq\emptyset$.
\end{enumerate}

It is well-known that the signed cocircuits satisfy the axioms for signed circuits and therefore define another matroid structure on the same set $\wt E$. This is called the \emph{dual matroid} and we denote it by $\wt M^\ast$. Important examples are the \emph{cographic matroids} which are the duals of graphic matroids.

We will drop the word ``signed'' and simply talk about circuits and cocircuits in $E$ and in $\wt E$. Bases and independent sets in $\wt E$ are simply subsets which map monomorphically onto bases and independent sets in $E$.

\begin{defn} If $(\wt E,\wt M)$ is a freely oriented matroid on a set $E$ and $E'$ is a subset of $E$, then the \emph{coinduced} freely oriented matroid on $E'$ is $(\wt E',\wt M')$ where
\begin{enumerate}
\item $\wt E'=p^{-1}(E')\subseteq \wt E$ and
\item $\wt M'=\wt M\cap\cP(\wt E')$ where $\cP(\wt E')$ is the power set of $\wt E'$. I.e., the circuits of $\wt E'$ are defined to be the circuits of $\wt E$ which lie in $\wt E'$.
\end{enumerate}
We also write
\[
    \wt M|E':=(\wt E',\wt M')
\]
and call this a \emph{submatroid} of $(\wt E,\wt M)$.
\end{defn}

Let $\cM(r,\ZZ)$ denote the category whose objects are freely oriented matroids of rank $r$ and whose morphisms are $\ZZ/2$ equivariant monomorphism
\[
    \wt\f:\wt E_1\to \wt E_2
\]
so that $C\subseteq \wt E_1$ is a circuit if and only if its image in $\wt E_2$ is a circuit. In other words $\wt\f$ is an isomorphism of $(\wt E_1,\wt M_1)$ with a submatroid of $(\wt E_2,\wt M_2)$.

\begin{defn} A \emph{representation} of $(\wt E,\wt M)$ over $\RR$ is a function
\[
    \r:\wt E\to V
\]
from $\wt E$ to a real vector space $V$ so that
\begin{enumerate}
\item $\r(e)+\r(\ov e)=0$ for all $e\in\wt E$
\item For any circuit $C$ we have a linear relation
\[
    \sum_{e_i\in C}a_i\r(e_i)=0
\]
where $a_i$ are all positive.
\item $\r$ maps any independent subset of $\wt E$ monomorphically onto a linearly independent subset of $V$.
\item The image of $\wt E$ spans $V$.
\end{enumerate}
The representation $\r$ is called \emph{unipotent} if the images of any two bases have the same volume in $V$. A unipotent representation of $(\wt E,\wt M)$ over $\ZZ$ is a function $\r:\wt E\to L$, where $L$ is a free abelian group, satisfying the above 4 conditions plus the additional condition that every basis maps to a set of free generators of $L$.
\end{defn}

\begin{defn}
A \emph{regular freely oriented matroid} (ROM) is defined to be a freely oriented matroid $(\wt E,\wt M)$ for which there exists a unipotent representation $\r:\wt E\to V$ over $\RR$.
\end{defn}

It is well-know and easy to prove that every submatroid of a ROM is a ROM.

Let $\cR(r,\ZZ)$ denote the full subcategory of $\cM(r,\ZZ)$ of regular freely oriented matroids. Then, in analogy with the binary case, we want to investigate the first integral homology functor:
\[
    H_1:\cR(r,\ZZ)\to \cG\el(r,\ZZ).
\]


\subsection{Homology of a freely oriented matroid}

\begin{defn}
The \emph{first homology group} $H_1(\wt M)=H_1(\wt E,\wt M)$ of a freely oriented matroid $(\wt{E},\wt M)$ on $E$ is defined to be the additive group generated by the elements of $\wt E$ modulo the following relations.
\begin{enumerate}
\item $e+\ov{e}=0$ for all $e\in \widetilde E$.
\item $\sum_{e\in C}e=0$ for all circuits $C\in \wt M$.
\end{enumerate}
Thus the elements of $H_1(\wt M)$ are equivalence classes $\left[\sum n_ie_i
\right]$ of additive linear combinations of elements of $\wt E$.
\end{defn}

If $\wt M'$ is a submatroid of $\wt M$, we get an induced map in homology:
\[
    H_1(\wt M')\to H_1(\wt M).
\]

\begin{prop}\begin{enumerate}\label{rank of H1}
\item If $\wt M'$ is a submatroid of $\wt M$ of full rank then the induced map in homology is an epimorphism.
\item $H_1(\wt M)$ is generated by any basis for $\wt E$
\item The rank of $H_1(\wt M)$ is bounded by the rank of the matroid $\wt M$.
\end{enumerate}
\end{prop}

\begin{proof} Since $\wt E'$ has full rank it contains a basis $B$ for $\wt E$. Take any $e\in \wt E\delete (B\cup \ov B)$. Then $B\cup\ov B\cup \{e\}$ contains a circuit $C$ which is not contained in $B\cup\ov B$. So, $e\in C$ and $C\delete e\subseteq B\cup \ov B$. Therefore, $[e]\in H_1(\wt M)$ is a linear combination of $[b]$, $b\in B$ which lie in the image of $H_1(\wt M')$. This proves (1) and (2). The rank statement follows.
\end{proof}

\begin{defn}
For any additive group $A$ we define the \emph{cohomology} group $H^1(\wt M;A)$ of $(\wt E,\wt M)$ with coefficients in $A$ to be the group of all functions $f:\widetilde E\to A$ satisfying the following conditions.
\begin{enumerate}
\item $f(\ov{e})=-f(e)$ for all $e\in \wt E$.
\item $\sum_{e\in C}f(e)=0$ for all circuits $C\in\wt M$.
\end{enumerate}
We call such functions $f$ \emph{cocycles} on $\wt E$ with coefficients in $A$.
\end{defn}

Of particular interest is the \emph{universal cocycle}
\[
    \iota: \wt E\to  H_1(\wt M)
\]
given by $\iota(e)=[e]$. This very clearly has the following universal property.

\begin{prop}
Given any cocycle $f:\wt E\to A$, there exists a unique homomorphism $\what f:H_1(\wt M)\to A$ so that the following diagram commutes.
\[
\xy
\xymatrix{
\wt E\ar[r]^{\iota\quad }\ar[dr]_f& H_1(\wt M)\ar[d]^{\what f}\\
 & A
}
\endxy
\]
Consequently, the correspondence $f\leftrightarrow\what f$ gives an isomorphism
\[
    H^1(\wt M;A)\cong \Hom(H_1(\wt M),A).
\]
\end{prop}

If $\wt M'$ is a submatroid of $\wt M$, we get an induced map in cohomology
\[
    H^1(\wt M;A)\to H^1(\wt M';A)
\]
for any abelian group $A$. We call this the \emph{restriction map} since it is given by restricting any cocycle $f:\wt E\to A$ to $\wt E'$.

The following statement follows immediately from the propositions above.

\begin{cor}
If $(\wt E',\wt M')$ is a submatroid of $(\wt E,\wt M)$ of full rank then the induced map in cohomology is a monomorphism for any $A$. I.e., any cocycle on $\wt E$ is uniquely determined by its restriction to $\wt E'$.
\end{cor}



\subsection{Regular matroids}

\begin{lem}\label{regular rep is a cocycle} Suppose that $(\wt E,\wt M)$ is a \emph{regular freely oriented matroid} on $E$. Then any unipotent representation $\r:\wt E\to V$ over $\RR$ is a cocycle.
\end{lem}

\begin{proof}
Take a circuit $C=\{e_0,\cdots,e_m\}$. Then, there exist positive real numbers $a_0,\cdots,a_m$ so that
\[
    \sum a_i\r(e_i)=0.
\] Since $C\delete e_0$ is independent, we can extend $e_1,e_2,\cdots,e_m$ to a basis $B$ by adding elements $b_1,\cdots,b_{r-m}$. Then the dual basis element $f:V\to \RR$ which sends $\r(e_1)$ to 1 and the image under $\r$ of all other basis elements to 0, when evaluated on the relation above gives:
\[
    a_1+a_0f(\r(e_0))=0.
\]
Since $\r$ is unipotent we have $f(\r(e_0))=0,1$ or $-1$. Since $a_0,a_1$ are positive, the only possibility is: $f(\r(e_0))=-1$ and $a_0=a_1$. Since the ordering of the elements of $C$ was chosen arbitrarily, we conclude that the coefficients $a_i$ must all be equal. This implies that $\r(C)=0$. Since this holds for every circuit $C$, we conclude that $\r$ is a cocycle.
\end{proof}

We have the following immediate consequence which is a special case of the well-known ``projective rigidity'' of regular matroids.

\begin{thm}\label{thm:H1(M)=L}
Given any unipotent representation $\r:\wt E\to V$ of a regular matroid on $E$, let $L$ be the additive subgroup of $V$ generated by the image of $\r$. Then there exists a unique isomorphism $\what\r:H_1(\wt M)\to L$ so that the following diagram commutes.\[
\xy
\xymatrix{
\wt E\ar[r]^{\iota\quad }\ar[dr]_\r& H_1(\wt M)\ar[d]^{\what\r}_\cong\\
 & L
}
\endxy
\]
\end{thm}

\begin{thm}
A freely oriented matroid $(\wt E,\wt M)$ of rank $r$ is regular if and only if $H_1(\wt M)\cong \ZZ^r$.
\end{thm}

\begin{proof}
If $\wt M$ is regular let $L$ be as given in the above theorem. Then $L$ is a torsion-free abelian group of rank at least $r$ since it spans $V\cong\RR^r$. But the rank of $H_1(\wt M)$ is at most $r$. So,
\[
    H_1(\wt M)\cong L\cong\ZZ^r.
\]

Conversely, suppose that $H_1(\wt M)\cong\ZZ^r$. Then we claim that the universal cocycle
\[
    \iota:\wt E\to V=H_1(\wt M)\otimes\RR
\]
is a unipotent representation of $\wt E$. The fact that $\iota$ is a representation follows from the fact that circuits go to zero and bases go to spanning sets by Proposition \ref{rank of H1}. Since any basis generates $H_1(\wt M)$, each basis must have volume equal to the covolume of the lattice $H_1(\wt M)$ in $V$. Therefore, $\iota$ is a unipotent representation.
\end{proof}

If we use cocircuits we get a symmetrical characterization of regular freely oriented matroids. For any two subsets $A,B$ of $\wt E$ the scalar product $\<A,B\>\in\ZZ$ is defined by
\[
    \<A,B\>=|A\cap B|-|A\cap\ov B|
\]
and the \emph{dual} of $A$ is defined to be the function $A^\vee:\wt E\to \ZZ$ given by
\[
    A^\vee(B)=\<A,B\>.
\]

\begin{thm} Let $(\wt M,\wt E)$ be a freely oriented matroid of rank $r$ on a set $E$ with $n$ elements. Then
the following are equivalent.
\begin{enumerate}
\item $H_1(\wt E)\cong\ZZ^r$ (i.e., $\wt M$ is regular).
\item The dual $C^\vee$ of any cocircuit $C$ is a cocycle.
\item For any circuit $C$ and cocircuit $C^\ast$, $\<C, C^\ast\>=0$.
\item $H_1(\wt E^\ast)\cong\ZZ^{n-r}$ (i.e., $\wt M^\ast$ is regular).
\end{enumerate}
\end{thm}

\begin{proof}
(2) and (3) are obviously equivalent. Also, (4) is dual to (1) and condition (3) is self-dual. So, it suffices to show that (1) is equivalent to (2).

$(2)\then(1)$ Let $B=\{b_1,\cdots,b_r\}$ be a basis in $\wt E$. Then for each $b_i\in B$ there is a cocircuit $C_i\subseteq \wt E$ so that $C_i\cap (B\cup\ov B)=b_i$. By (2), the dual $C_i^\vee$ of $C_i$ gives a cocycle
\[
    C_i^\vee: H_1(\wt E)\to\ZZ
\]
so that $C_i^\vee(b_j)=\delta_{ij}$. Thus $(C_1^\vee,\cdots,C_r^\vee)$ is a left inverse for the surjection $\ZZ^B\onto H_1(\wt E)$ making that map an isomorphism.

$(1)\then(2)$ The universal cocycle
\[
    \iota:\wt E\to H_1(\wt M)\hookrightarrow V=H_1(\wt M)\otimes\RR
\]
is a representation of $\wt E$ which sends every basis to a generating set for the lattice $H_1(\wt M)$ and therefore has volume 1. This implies that the projection
\[
    \pi:V\to\RR
\]
along any hyperplane $H=\ker\pi$ spanned by elements in the image of $\wt E$
has values $0,1,-1$ (after normalizing). Therefore, this projection is equal to the dual of the corresponding cocircuit $C^\ast=\iota^{-1}\pi^{-1}(1)$ making it a cocycle. Furthermore, all cocircuits $C^\ast$ occur in this way since the set of all $e\in \wt E$ which are disjoint from $C^\ast$ and $\ov C^\ast$ span a hyperplane.
\end{proof}


\subsection{The fiber poset for ROMs}

Homology is a functor from $\cR(r,\ZZ)$ the category of all ROMs of rank $r$ to the category of free abelian groups of rank $r$ and isomorphisms:
\[
    H_1:\cR(r,\ZZ)\to \cG\el(r,\ZZ).
\]
In complete analogy with Definition \ref{def:fiber poset for regular binary matroids} and Proposition \ref{prop: fiber poset for regular matroids} we have the following.

\begin{defn} Let $IR(r,\ZZ)$ denote the poset of all subsets $\wt E$ of $\ZZ^r$ which are images of unitary representations of regular freely oriented matroid of rank $r$. Equivalently, $IR(r,\ZZ)$ is the poset of all subsets $\wt E$ of $\ZZ^r$ ordered by inclusion, which have the following properties.
\begin{enumerate}
\item $\wt E$ does not contain $0$
\item $\wt E=-\wt E$
\item $\wt E$ generates $\ZZ^r$ additively.
\item Any set of $r$ linearly independent elements of $\wt E$ generates $\ZZ^r$.
\end{enumerate}
\end{defn}

As in Definition \ref{mod2 homologically marking}, a \emph{homological marking} on a regular oriented matroid is an isomorphism $H_1(\wt M)\cong\ZZ^r$ or, equivalently, a unipotent representation $\wt E\to\ZZ^r$. The category of homologically marked freely oriented regular matroids of rank $r$ is the fiber category of the homology functor and is denoted $\cI\cR(r,\ZZ)$

\begin{prop} The geometric realization of $IR(r,\ZZ)$ is the homotopy fiber of the mapping\[
  B\cR(r,\ZZ)\to BGL(r,\ZZ)
\]
induced by the first homology functor $H_1$.
Equivalently, $|IR(r,\ZZ)|$ is homotopy equivalent to the covering space of $B\cR(r,\ZZ)$ classified by the functor $H_1$. Equivalently, $|IR(r,\ZZ)|\simeq B\cI\cR(r,\ZZ)$.
\end{prop}

\begin{lem}\label{lem:reduction of IR mod 2}
Reduction modulo $2$ gives an order preserving functor $p:IR(r,\ZZ)\to IR(r,\FF_2)$. Furthermore, each $\wt E$ in $IR(r,\ZZ)$ is a two fold covering of $E=p(\wt E)$.
\end{lem}

\begin{proof}
Reduction mod 2 is obviously an order preserving map. It is onto by the definition of a regular matroid. To prove the last statement, suppose that $\wt E$ has two elements $e_1,e_2$ which are congruent modulo 2. Then they cannot be part of any basis for $\ZZ^r$. Therefore, they must be linearly dependent (since any independent set is contained in a basis). So, $ae_1=be_2$. If we take any basis containing $e_1$ and replace $e_1$ by $e_2$ then the determinant of the matrix will change by a factor of $a/b$. Therefore, $a/b=\pm1$ and $e_1=\pm e_2$ are required.
\end{proof}

We now look at the easy case $r=2$ and determine the homotopy type of the poset $IR(2,\ZZ)$ and category $\cR(2,\ZZ)$ of regular matroids of rank 2.

\begin{thm}
$|IR(2,\ZZ)|$ is contractible and therefore $
B\cR(2,\ZZ)\simeq BGL(2,\ZZ).
$
\end{thm}

\begin{proof} We use the lemma above. The poset $IR(2,\FF_2)$ has only four elements. They are the three bases for $\FF_2^2$ each with 2 elements and the single circuit with three elements. Thus each $\wt E\in IR(2,\ZZ)$ is either minimal with 4 elements or maximal with 6 elements:
\begin{enumerate}
\item If $\wt E$ has 4 elements they must be $\pm a,\pm b\in \ZZ^2$ where $a,b$ forms a basis for $\ZZ^2$.
\item If $\wt E$ has 6 elements they must be $\pm a,\pm b,\pm c\in\ZZ^2$ where $a+b+c=0$, i.e. $a,b,c$ form a triangle, and this triangle must have area $1/2$.
\end{enumerate}

In case (2), we claim that the triangle is either
\begin{enumerate}
\item[(2a)] a right triangle with sides of length $1,1,\sqrt2$ (only $\{\pm(1,0),\pm(0,1),\pm(1,1)\}$ and $\{\pm(1,0),\pm(0,1),\pm(1,-1)\}$ are of this kind)
\item[(2b)] a triangle with one obtuse angle.
\end{enumerate}
The proof of this is very simple. Take a triangle with no obtuse angles and let $s$ be the length of the shortest side. If $s=1$ we must be in case (2a). Otherwise $s\ge\sqrt2$ and a simple drawing shows that the smallest possible triangle with these properties is an equilateral triangle with area
\[
    area\ \ge \frac{\sqrt3}4 s^2\ge \frac{\sqrt3}2>\frac12
\]
which is a contradiction.

Every matroid of type 2 contains exactly three submatroids of type 1 given by each of the 3 bases of $E$. Each matroid of type 1 is a submatroid of exactly 2 matroids of type 2 given by adding either $\pm(a+b)$ or $\pm(a-b)$ if $a,b$ form a basis for $\wt E$. Thus, $|IR(2,\ZZ)|$ is a graph having 2-valent and 3-valent vertices.
\begin{figure}[htbp]
\begin{center}
%
{
\setlength{\unitlength}{1.5cm}
{\mbox{
\begin{picture}(4,2)
      \thicklines
  %
    \put(2,1){\circle*{.1}} 
    \put(1.45,.75){$\{\pm a,\pm b\}$}
    \put(1.3,1){\circle*{.1}} 
    \put(1.3,1){\line(1,0){1.5}
    }
    \put(2.8,1){\circle*{.1} } 
    \put(3.1,.95){$\{\pm a,\pm b, \pm (a+b)\}$}
    \put(-1.3,.95){$\{\pm a,\pm b, \pm (a-b)\}$}
    \put(2.7,1){
    \line(1,1){.8}}
    \put(3.57,1.8){\circle*{.1}} 
    \put(3.57,1.8){\line(0,1){.4}} 
    \put(3.57,1.8){\line(1,0){.5}} 
    \put(3.57,0.2){\circle*{.1}} 
    \put(3.57,0.2){\line(0,-1){.3}} 
    \put(3.57,0.2){\line(1,0){.5}}  
   \put(3.18,1.4){\circle*{.1}} 
   \put(3.18,.6){\circle*{.1}} 
    \put(.9,.6){\circle*{.1}} 
   \put(.9,1.4){\circle*{.1}} 
    \put(.5,1.8){\circle*{.1}} 
    \put(.5,1.8){\line(0,1){.4}}  
    \put(.5,1.8){\line(-1,0){.5}}  
    \put(.5,.2){\circle*{.1}} 
    \put(.5,.2){\line(0,-1){.3}} 
    \put(.5,.2){\line(-1,0){.5}} 
   \put(2.7,1){
     \line(1,-1){.8}   }
    \put(1.23,1){
    \line(-1,1){.8}}
   \put(1.23,1){
     \line(-1,-1){.8}   }
\end{picture}}
}}
\label{default}
\end{center}
\end{figure}
We claim that this is an infinite tree and thus contractible. To prove this, we define a Morse function $f$ on the graph by taking each matroid of type 2 (a 3-valent vertex of the graph) to the length of its maximal element. At every 3-valent vertex of type (2b) of the graph, the function decreases in only one direction, namely the one in which the side $a+b$ opposite the obtuse angle has been exchanged with $a-b$. The unique  pair of vertices without this property is given by the two right triangles given in Case (2a) above. Thus, the space $|IR(2,\ZZ)|$ collapses to the line segment containing these two vertices (and the intermediate 2-valent vertex representing the matroid $\{\pm(1,0),\pm(0,1)\}$) and is thus contractible.
\end{proof}

The action of $GL(2,\ZZ)$ on this infinite tree is well-known and used in many application \cite{SerreTrees}. For example, the level $2$ congruence subgroup
\[
    \Gamma(2):=\ker(GL(2,\ZZ)\to GL(2,\FF_2))
\]
acts on this graph and the quotient space is $|IR(2,\FF_2)|$ which is a tree with $4$ vertices, the 4 elements of the poset $IR(2,\FF_2)$. By Lemma \ref{lem:reduction of IR mod 2} the action of $\Gamma(2)$ on $|IR(2,\ZZ)|$ is free and transitive on the 2-valent vertices and has three orbits of 2-valent vertices, each with 2-element stabilizer subgroups. By van Kampen's theorem this implies that $\Gamma(2)$ is a triple free product:
\[
    \Gamma(2)\cong\ZZ/2\ast \ZZ/2\ast\ZZ/2.
\]


%
%

\section{Graphs}

Graphic and cographic matroids are the main examples of regular matroids. By definition, these matroids make up the image of functors from two categories of graphs into the category $\cR(r,\ZZ)$. We want to examine these two graph categories and the corresponding matroid categories.


\subsection{Graph monomorphisms}

The first category is the category of finite connected graphs and graph monomorphisms. 

\begin{defn}
By a (finite) \emph{graph} we mean a finite set of vertices $G_0$, a finite set of edges $G_1$, a two-fold covering $\wt G_1\to G_1$ called the set of \emph{half-edges} and an \emph{incidence map} $\d:\wt G_1\to G_0$.
\end{defn}

One example is the \emph{complete graph} $K(S)$. For any set $S$, this is the graph whose vertex set is $S$, edge set $E$ is the set of all subsets of $S$ with exactly two elements and half-edges are ordered pairs of two distinct elements of $S$.

For any graph $G$, the set $G_0\coprod \wt G_1$ has two endomorphism $\t$ and $\d$. On $\wt G_1$, $\t$ is the free involution which switches the two sheets over $G_1$ and $\d$ is the incidence map $\wt G_1\to G_0$. On $G_0$ both maps are defined to be the identity map ($\t(v)=v=\d v$ for all $v\in G_0$). A general morphism of graphs $G=(G_0,\wt G_1)\to G'=(G_0',\wt G_1')$ is any mapping
\[
    G_0\smallcoprod \wt G_1\to G_0'\smallcoprod \wt G_1'
\]
which commutes with both $\t$ and $\d$.

\begin{defn}
A \emph{graph monomorphism} is a morphism of graphs
\[
    f:G\to G'
\]
which is an isomorphism of $G$ onto a spanning subgraph of $G'$. Thus $f$ is a monomorphism of sets
\[
    f:G_0\smallcoprod\wt G_1\to G_0'\smallcoprod\wt G_1'
\]
so that $f(G_0)=G_0'$ and $f(\wt G_1)\subseteq \wt G_1'$.
\end{defn}

Recall that a \emph{spanning} subgraph is a subgraph having all the vertices of the original graph. Every connected graph contains a spanning tree as a subgraph.

\begin{prop} Every graph monomorphism between connected graphs is a composition of graph isomorphisms and inclusion maps
\[
    G\del e\hookrightarrow G.
\]
where $G\del e$ is $G$ with a nonseparating edge $e$ deleted.
\end{prop}

Suppose that $G=(G_0,\wt G_1)$ and $G'=(G_0,\wt G_1')$ are two graphs with the same set of vertices. Then we define the union $G\cup G'$ to be the graph $(G_0,\wt G_1\coprod \wt G_1')$. This graph contains both $G$ and $G'$ as subgraphs. Strictly speaking, $G\cup G'$ contains two subgraphs isomorphic to $G,G'$ (e.g., if $G=G'$). I.e., we have graph monomorphisms
\[
    G\to G\cup G'\ot G'.
\]

\begin{defn} Let $m\cG(r)$ denote the category whose objects are connected graphs with $r+1$ vertices and whose morphisms are graph monomorphisms. Let $\cS(r+1)$ denote the category whose objects are sets with cardinality $r+1$ and whose morphisms are bijections. Then we define the \emph{vertex functor}
\[
    V:m\cG(r)\to\cS(r+1)
\]
to be the functor which sends $G=(G_0,\wt G_1)$ to the vertex set $G_0$.
\end{defn}

\begin{thm}\label{thm:mG(r)=Sr+1}
The vertex functor gives a homotopy equivalence of classifying spaces
\[
    Bm\cG(r)\simeq B\cS(r+1)\simeq BS_{r+1}
\]
where $S_n$ is the symmetric group on $n$ letters.
\end{thm}

\begin{proof}
A right inverse to the vertex functor is given by the complete graph functor $K:\cS(r+1)\to m\cG(r)$ which sends each set to the complete graph on that set. Since the vertex set of the complete graph is the original set, we have $V\circ K=id$. Conversely, $K\circ V\simeq id$ by the following construction:
\[
    KV(G)\to KV(G)\cup G\ot G.
\]
\end{proof}

We have a functor
\[
    M:m\cG(r)\to \cR(r,\ZZ)
\]
which sends $G=(G_0,\wt G_1)$ to $\wt G_1$ considered as a graphic matroid as explained earlier. The composition of functors
\[
    \cS(r+1)\xrarrow{K} m\cG(r)\xrarrow{M} \cR(r,\ZZ)\xrarrow{H_1} \cG\el(r,\ZZ)
\]
sends each set $S$ to the kernel of the homomorphism
\[
    \ZZ^S\to \ZZ
\]
sending $x\in\ZZ^S$ to $\sum x_s$, the sum of the coordinates of $x$.

Let $\cI m\cG(r)$ be the category of homologically marked graphs in $m\cG(r)$ where by a \emph{homological marking} on $G$ we mean a homological marking on the corresponding graphic matroid. Thus $\cI m\cG(r)$ is the pull-back in the following diagram of categories.
\[
\xy
\xymatrix{
\cI m\cG(r)\ar[r]\ar[d] &\cI \cR(r,\ZZ)\ar[d]\\
m\cG(r)\ar[r] &\cR(r,\ZZ)
}
\endxy
\]

\begin{cor}
The classifying space of $\cI m\cG(r)$ is homotopy equivalent to the discrete set of cosets of $W=S_{r+1}$ in $GL(r,\ZZ)$:
\[
    B\cI m\cG(r)\simeq GL(r,\ZZ)/S_{r+1}.
\]
\end{cor}

Results about the category of graphs do not translate easily into results about graphic matroids because of the following observation.

\begin{prop}
The functors $M:m\cG(r)\to \cR(r,\ZZ)$ and $M:\cI m\cG(r)\to \cI\cR(r,\ZZ)$ are faithful for all $r$ but not full for $r\ge2$.
\end{prop}

\begin{proof}
The functor $M$ is faithful because any graph monomorphism of connected graphs is determined by its restriction to the half-edge set. For the second part, consider the involution on $M(G)$ which reverses the sign of every element. This is not induced by a graph monomorphism if $r\ge2$.
\end{proof}


\subsection{Graph epimorphisms}

By a {graph epimorphism} we mean a morphism of graphs $f:G\to G'$ which is a composition of isomorphisms and collapsing maps
\[
    G\to G/e
\]
where $G/e$ is $G$ with the edge $e$ collapsed where $e$ is not a loop. This is the same as saying that the inverse image of every edge of $G'$ is an edge of $G$ and the inverse image of every vertex of $G'$ is a subgraph of $G$ which is a tree.

\begin{defn}
A \emph{graph epimorphism} is a morphism of graph $f:G\to G'$ which is given by an epimorphism
\[
    f:G_0\smallcoprod\wt G_1\to G_0'\smallcoprod\wt G_1'
\]
satisfying the following.
\begin{enumerate}
\item $f(G_0)=G_0'$.
\item Every half-edge of $G'$ is the image of exactly one half-edge in $G$.
\item The inverse image of every vertex in $G'$ is a tree in $G$.
\end{enumerate}
\end{defn}

Let $e\cG(r)$ be the category of connected graphs with Euler characteristic $1-r$ (this is the number of vertices minus the number of edges) with graph epimorphisms as morphisms. Let
\[
    M^\ast:e\cG(r)\to \cR(r,\ZZ)
\]
be the contravariant functor which associates the cographic matroid $M^\ast(G)$ to $G$. (This is the dual of the graphic matroid $M(G)$ of $G$. The circuits are the minimal separating sets of edges.) A graph epimorphism $f:G\to G'$ induces a morphism of freely oriented matroids $M^\ast(G')\to M^\ast(G)$ by sending $e\in\wt E'$ to $f^{-1}(e)\in\wt E$.

We have the following well-known theorem. (See \cite{Igusa:book}, p.307.)

\begin{thm}\label{thm:eG(r)=Out(Fr)}
$Be\cG(r)\simeq BOut(F_r)$, the classifying space of the group $Out(F_r)$ of outer automorphisms of the free group $F_r$ on $r$ letters.
\end{thm}

Since the abelianization $H_1(Out(F_r))$ of $Out(F_r)$ is isomorphic to $ GL(r,\ZZ)$, we get the following corollary.

\begin{cor}
The functor $H_1:\cR(r,\ZZ)\to \cG\el(r,\ZZ)$ induces an epimorphism on fundamental groups:
\[
    \pi_1 \cR(r,\ZZ)\onto GL(r,\ZZ).
\]
\end{cor}

S\o ren Galatius \cite{Galatius:Free} has recently computed the homology of this category in the stable range (in degree $\le r/2-2$) and shown that it has the same homology as the symmetric group. By Theorems \ref {thm:mG(r)=Sr+1} and \ref{thm:eG(r)=Out(Fr)} this implies the following.

\begin{thm}\label{thm:graphs and cographs are homologically equivalent}
$e\cG(r)$ has the same homology as $m\cG(r)$ in a stable range ($k\le r/2-2$):
\[
    H_k(Be\cG(r))\cong H_k(Bm\cG(r))
\]
\end{thm}

This leads to the following question:

\begin{quest}
Do $M:m\cG(r)\to\cR(r,\ZZ)$ and $M^\ast:e\cG(r)\to\cR(r,\ZZ)$ induce isomorphisms in homology in a stable range?
\end{quest}

Let $\cI e\cG(r)$ be the category of cohomologically marked connected graphs with Euler characteristic $1-r$ where, by a \emph{cohomological marking} we mean a marking of the cographic matroid $M^\ast(G)$. In standard terminology, this is an ordered basis for the first singular cohomology group
\[
    H^1(|G|;\ZZ)\cong \ZZ^r
\]
of the topological realization $|G|$ of the graph $G$. Then $GL(r,\ZZ)$ acts freely on $\cI e\cG(r)$ and the quotient category is $e\cG(r)$. So, we get a homotopy fiber sequence of classifying spaces:
\[
	B\cI e\cG(r)\to B e\cG(r)\to B\cG\el(r,\ZZ)
\]
Which, by Theorem \ref{thm:eG(r)=Out(Fr)} above, implies the following.

\begin{cor}\label{cor:IOut is IeG}
$B\cI e\cG(r)$ is homotopy equivalent to the classifying space of the kernel $IOut(F_n)$ of the epimorphism $Out(F_n)\to GL(n,\ZZ)$.
\end{cor}

As we will explain in the next section, we are very interested in the cohomology of $IOut(F_n)$ which, by this remark, is the same as the cohomology the category $\cI e\cG(r)$. My idea is that the corresponding category of regular matroids $\cI \cR(r,\ZZ)$ may be easier to study. As in the graphic matroid case, the following observation shows that category of cographic matroids is not the same as $e\cG(r)$.

\begin{prop}
The functors $M^\ast:e\cG(r)\to \cR(r,\ZZ)$ and $M^\ast:\cI e\cG(r)\to \cI\cR(r,\ZZ)$ are faithful but not full for $r\ge1$.
\end{prop}

\begin{proof}
As before, these functors are faithful since an epimorphism of connected graphs is uniquely determined by the operation which takes each half-edge in the target to its inverse image. The functors are not full since the automorphism of $M^\ast(G)$ given by reversing the sign of each element is not, in general, induced by an automorphism of $G$.
\end{proof}

%
%

\section{Relation to higher Reidemeister torsion}

The main motivation behind this investigation of the homotopy type of the category of regular matroids is the speculation that the higher Reidemeister torsion invariants can be defined on the category $\cI\cR(r,\ZZ)$ where it might be easier to calculate and easier to generalize to the complexified case. In this section I will review what is known about higher torsion and how I expect this to extend to the category of matroids. I will use the axiomatic approach to higher torsion developed in \cite{IAxioms0}.

\subsection{Review of axiomatic higher torsion}

We consider smooth $(C^\infty)$ manifold bundles $p:E\to B$ where the fiber $X$ is an oriented closed manifold and $B$ is connected. We call such a bundle \emph{unipotent} if it satisfies the property that as a $\pi_1B$-module, the rational homology of the fiber $H_m(X;\QQ)$ in every degree $m$ has a filtration by submodules so that the action of $\pi_1B$ on each subquotient is trivial. For example, if $B$ is simply connected this condition is automatically satisfied.

By a \emph{real characteristic class} for unipotent bundles we mean a cohomology class $\t(E)\in H^\ast(B;\RR)$ with the following naturality property. Suppose that $f:B'\to B$ is a smooth mapping of smooth manifolds. Then we have the pull-back bundle $f^\ast E\to B'$ given by $f^\ast E =\{(e,b)\in E\times B'\st p(e)=f(b)\}$. Then we require that the characteristic class for $E$ pulls back to that of $f^\ast E$:
\[
	f^\ast(\t(E))=\t(f^\ast E)\in H^\ast(B';\RR)
\]

\begin{defn} For a fixed positive integer $k$ we define a degree $4k$ \emph{higher torsion invariant} to be any real characteristic class $\t$ for unipotent smooth bundles
\[
	\t(E)\in H^{4k}(B;\RR)
\]
which satisfies the following two conditions.
\begin{enumerate}
\item (\emph{additivity}) Suppose that $E$ is the union of two smooth manifold bundles $E_i\to B$ whose fibers $X_i$ are compact manifolds with the same boundary $\d X_1=\d X_2$. For each $i$, let $DE_i\to B$ be the fiberwise double of $E_i$. Then
\[
	{\t(E)=\frac12\t(DE_1)+\frac12\t(DE_2)}
\]
\item (\emph{transfer}) Suppose that $D\to E$ is an oriented linear $m$-sphere bundle and let $\t_E(D)\in H^{4k}(E;\RR), \t_B(D)\in H^{4k}(B;\RR)$ be the higher torsion of $D$ as a bundle over $E$ and $B$ respectively. Then
\[
	{\t_B(D)=\chi(S^m)\t(E)+tr_B^E(\t_E(D))}
\]
where the transfer $tr_B^E:H^\ast(E)\to H^\ast(B)$ is defined below.
\end{enumerate}
\end{defn}

Suppose that $E\to B$ is a smooth bundle with fiber $X$ a closed oriented manifold of dimension $n$. Then the relative tangent bundle $TX$ is the tangent bundle of all of the fibers considered as a bundle over $E$. This is also the kernel of the derivative map $TE\to TB$ induced by $p:E\to B$. Let $e(TX)\in H^n(E)$ be the Euler class of this $n$ dimensional vector bundle. Then the transfer $tr^E_B:H^\ast (E)\to H^\ast(B)$ is defined by
\[
	tr^E_B(x)=p_\ast(x\cup e(TX))
\]
where $p_\ast:H^{m+n}(E)\to H^m(B)$ is given in deRham cohomology by ``integration along the fibers'' and in integer cohomology by taking the cohomology Serre spectral sequence of the bundle $E\to B$ and observing that $H^{m+n}(E)$ maps onto $E_\infty^{m,n}$ which is a subgroup of $H_2^{m,n}=H^m(B;H^n(X))\cong H^m(B)$. 

Any higher torsion invariant can be extended to unipotent smooth bundles $E\to B$ whose fibers $X$ are compact manifolds with boundary by the formula
\[
	\t(E)=\frac12\t(DE)+\frac12\t(\d^v E)
\] 
where $DE$ is the fiberwise double of $E$ and $\d^v E$ denotes the subbundle of $E$ whose fiber is $\d X$, the boundary of $X$.

\begin{defn}
A higher torsion invariant $\t$ is called \emph{exotic} if it has the property that
\[
	\t(D)=\t(E)
\]
for any linear disk bundle $D$ over $E$.
\end{defn}

It was shown in \cite{I:ComplexTorsion} that the Igusa-Klein higher torsion invariant $\t^{IK}$ constructed using Morse theory in \cite{Igusa:book} (generalizing the construction and calculation in \cite{IK1Borel2}) is an exotic higher torsion invariant. One of the main results of the paper \cite{IAxioms0} is the following.

\begin{thm}
Any degree $4k$ exotic higher torsion invariant is a scalar multiple of $\t^{IK}$ for each $k>0$.
\end{thm}

In \cite{DWW}, another higher torsion invariant was constructed. Recently, these invariants were shown to be essentially equal.

\begin{thm}\cite{BDKW}
The Dwyer-Weiss-Williams smooth torsion satisfies the axioms for an exotic higher torsion invariant and is therefore equal to a scalar multiple of Igusa-Klein torsion.
\end{thm}

\subsection{Higher torsion on graph bundles}

A key observation about the higher torsion is its relationship to the Miller-Morita-Mumford class.

\begin{defn}
Suppose that $E\to B$ is a smooth bundle with closed oriented even dimensional fibers. Then the \emph{generalized Miller-Morita-Mumford class} $M_{2k}(E)\in H^{4k}(B;\RR)$ is defined by
\[
	M_{2k}(E)=\frac12tr^E_B((2k)!ch_{4k}(TX\otimes \CC))
\]
\end{defn}

When the fibers are surfaces, these are the standard Miller-Morita-Mumford classes:
\[
	M_{2k}(E)=\k_{2k}=p_\ast(e(T\Sig)^{2k+1})
\]
where $T\Sig=TX$ is the relative tangent bundle of $E$. These are also know as the (even) \emph{tautological classes}. It was conjectured by Klein and shown by Hain, Penner and myself that these are scalar multiples of the higher torsion invariants. In \cite{Igusa:book} the proportionality constant was computed:

\begin{thm}
Suppose that $E\to B$ is an oriented unipotent surface bundle. Then
\[
	\t^{IK}_{2k}(E)=\frac12(-1)^k\z(2k+1)\frac{\k_{2k}(E)}{(2k)!}\in H^{4k}(B;\RR)
\]
where $\z(s)=\sum \frac1{n^s}$ is the Riemann zeta function.
\end{thm}

Recall that the \emph{Torelli group} $I_g$ is the subgroup of the mapping class group of genus $g$ surfaces $\Sig_g$ consisting of all isotopy classes of automorphisms of $\Sig_g$ which induce the identity on $H_1(\Sig_g)$. The tautological bundle over $BI_g$ with fiber $\Sig_g$ is unipotent and therefore the above theorem applies with $B=BI_g$.

To pass to graphs and matroids we need surfaces with at least one marked point. The reason is that a punctured surface has a canonical embedded graph known as a \emph{fat graph}. (See \cite{Penner87}.) The mapping which sends a marked surface to a punctured surface and then to the graph gives a homomorphism of groups
\[
	\f: I_{g,s}\to IOut(F_{2g+s-1})
\]
where $I_{g,s}$ is the group of isotopy classes of automorphisms of $\Sig_g$ which fix $s$ points.

It was shown in \cite{IKW} and in \cite{Igusa:book} that higher torsion invariants $\t^{DWW}$ and $\t^{IK}$ respectively are defined on $IOut(F_n)$
\[
	\t_{2k}(IOut(F_n))\in H^{4k}(IOut(F_n);\RR)
\]
In \cite{Igusa:book} it was shown that the homomorphism $\f$ above respects these higher torsion classes in the following way. (Recall that, by \cite{BDKW}, $\t^{DWW}$ and $\t^{IK}$ are proportional.)

\begin{thm}[\cite{Igusa:book}] The homomorphism $\f: I_{g,s}\to IOut(F_{2g+s-1})$ sends $\t_{2k}$ to
\[
	\f^\ast\t_{2k}^{IK}(IOut(F_{2g+s-1}))=\frac{(-1)^k}{2(2k)!}\zeta(2k+1)\left(
	\k_{2k}-\sum_{i=1}^s2 e(\g_i)^{2k}
	\right)
\]
where $\g_i$ is the relative tangent bundle along the $i$-th marked point.
\end{thm}

This is Theorem 8.5.17 in \cite{Igusa:book} where the case $s=1$ is stated but the proof works for any $s$.

The odd tautological classes $\k_{2k+1}$ vanish on the Torelli group. However, it is unknown whether the even ones are trivial or not. The fact that they come from the cohomology of $\cI e\cG(2g+s-1)$ (whose classifying space is homotopy equivalent to $BIOut(F_{2g+s-1})$ by Corollary \ref{cor:IOut is IeG}) suggests that we can investigate this question on that category of graphs. For simplicity we now assume $s=1$.

We can go from the Torelli group to $IOut(F_{2k})$ to the corresponding category of regular matroids and obtain maps of classifying spaces:
\[
	BI_{g,1}\to BIOut(F_{2g})\simeq B\cI e\cG(2g)\to B\cI  \cR(2g,\ZZ)
\]
giving maps in cohomology
\[
	H^{4k}(I_{g,1};\RR)\ot H^{4k}(IOut(F_{2g});\RR)\ot H^{4k}(\cI \cR(2g,\ZZ);\RR)
\]
My hope is that the higher torsion invariants can be defined on the category $\cI\cR(n,\ZZ)$ in such a way that they pull back to the higher torsion invariants for $IOut(F_n)$. The lowest of these classes should live in degree $4$. Furthermore, I believe that a complexified version of this construction should give higher torsion invariants in all even degrees on some category of the form ``$\cI\cR(n,\ZZ[i])$'' or with $\ZZ[i]$ replaced by other rings with complex embeddings. This hypothetical element should be easier to find since the first case should be in degree 2.

%
%

\section{Proof of Theorem \ref{thm:IR(F2) is simply connected}}\label{last section}

In this section we show that the poset $IR(r,\FF_2)$ of all regular submatroids of $\FF_2^r$ of rank $r$ is simply connected for all $r\ge1$. The proof is in several steps. First we show that this poset is connected by showing that any regular matroid is connected by a ``regular geodesic'' to any basis for $\FF_2^r$. More generally any two regular matroids $E_1,E_2$ whose intersection $X=E_1\cap E_2$ is independent are connected by a regular geodesic. Next we show that any two regular geodesics are homotopic in $IR(r,\FF_2)$. Finally we show that any path is homotopic to a regular geodesic if there is a regular geodesic. The theorem follows easily from this even though there are examples (given below) where no regular geodesic exists.


\subsection{Existence of geodesics}

Suppose that $E_1,E_2$ are elements of the poset $IR(r,\FF_2)$. Then we define the \emph{(binary) distance} between $E_1,E_2$ to be:
\[
    d(E_1,E_2)=|E_1\backslash E_2|+|E_2\backslash E_1|.
\]
This is the length of the shortest path from $E_1$ to $E_2$ in the poset of all subsets of $\FF_2^r$. It is equal to the number of elements which must be deleted from $E_1$ plus the number of elements which must be added. If $d(E_1,E_2)=d$, we define a \emph{regular geodesic} from $E_1$ to $E_2$ to be a path $\g$ of length $d$ from $E_1$ to $E_2$ in $IR(r,\FF_2)$. This is a sequence of elements $\g(0),\g(1),\cdots,\g(d)\in IR(r,\FF_2)$ so that $\g(0)=E_1,\g(d)=E_2$ and $d(\g(i-1),\g(i))=1$ for $i=1,\cdots,d$.

Here is an example of regular matroids $E_1,E_2\in IR(r,\FF_2)$ which are not connected by a regular geodesic. Take the following binary matroid of rank $5$.
\[
    \left[
    \begin{array}{ccccc|cccc}
    1 & 0 & 0 & 0 & 0 & 0 & 0 &0 & 1 \\
    0 & 1 & 0 & 0 & 0 & 0 & 0 & 0 & 1 \\
    0 & 0 & 1 & 0 & 0 & 1 & 1 & 0 & 1\\
    0 & 0 & 0 & 1 & 0 & 1 & 0 & 1 & 1\\
    0 & 0 & 0 & 0 & 1 & 0 & 1 & 1 & 1
    \end{array}
    \right]
\]
 This matroid is not regular since it contains $F_7$ as a minor. However, if we delete either of the first two columns it will be regular since we can then clear the last column by pivoting. Let $E_1,E_2$ be the regular matroids obtained by deleting the first and second columns of this matroid respectively. Then $d(E_1,E_2)=2$ and there are exactly two paths of length 2 from $E_1$ to $E_2$ in the poset of all subsets of $\FF_2^5$. These paths pass through either $E_1\cup E_2$ which is not regular or $E_1\cap E_2$ which has rank 4. Therefore, there is no regular geodesic from $E_1$ to $E_2$ in this case. We need some hypotheses to guarantee the existence of regular geodesics.
 
 \begin{lem}\label{existence lemma}
 If $E_1,E_2\in IR(r,\FF_2)$ and $X=E_1\cap E_2$ is independent then there exists a regular geodesic from $E_1$ to $E_2$.
 \end{lem}
 
 This lemma will imply in particular that any basis for $F_2^r$ is connected to any point in $IR(r,\FF_2)$ by a regular geodesic. Therefore, $IR(r,\FF_2)$ is connected.
 
 \begin{proof}
Since $X$ is independent, it is contained in bases $B_1\subseteq E_1, B_2\subseteq E_2$. A regular geodesic from $B_1$ to $B_2$ is given by alternating between adding elements of $B_2\backslash B_1$ one at a time and deleting elements of $B_1\backslash B_2$ one at a time. By the basis exchange property of matroids, this can be done while keeping the rank equal to $r$. The matroids in this path are all regular since they have $r$ or $r+1$ elements. Together with geodesics $E_1\to B_1$, $B_2\to E_2$ this gives a regular geodesic from $E_1$ to $E_2$.
 \end{proof}
 

\subsection{Uniqueness of geodesics up to homotopy}

 \begin{lem}
 Any two regular geodesics from $E_1$ to $E_2$ are homotopic in $IR(r,\FF_2)$ (fixing the endpoints) assuming that $X=E_1\cap E_2$ is independent.
 \end{lem}
 
 \begin{proof}
 Let $|X|=n$ and $|E_i|=r+k_i$. Then $d(E_1,E_2)=2(r-n)+k_1+k_2$. We will prove the lemma by induction on the pair
 $
     (r-n,k_1+k_2).
 $
 
Suppose first that $n=r$. Then $X\in IR(r,\FF_2)$ and all regular geodesics from $E_1$ to $E_2$ lie inside the contractible set
\[
    C(X):=\{E\in IR(r,\FF_2)\st X\subseteq E\}
\]
So, the lemma holds.

Suppose $n<r$ and $k_1+k_2=0$. In other words, $E_1$ and $E_2$ are both bases. Then any geodesic from $E_1$ to $E_2$ starts by adding one element of $E_2\backslash E_1$ to $E_1$. If $e\in E_2\backslash E_1$ then $(E_1+e)\cap E_2=X+e$ is independent since it is a subset of $E_2$. Therefore, by induction on $r-n$, any two regular geodesics from $E_1+e$ to $E_2$ are homotopic. In other words, the homotopy class of any regular geodesic from $E_1$ to $E_2$ depends only on its first step. Let $\g(+e)$ denote any regular geodesic through $E_1+e$. If $e_1,e_2$ are two elements of $E_2\backslash E_1$ then $E_1+e_1+e_2$ is regular since it has only $r+2$ elements. The two paths $E_1\to E_1+e_1\to E_1+e_1+e_2$ and  $E_1\to E_1+e_2\to E_1+e_1+e_2$ are homotopic, since they lie in $C(E_1)$. Compose this with any regular geodesic from $E_1+e_1+e_2$ to $E_2$ to get a homotopy $\g(+e_1)\simeq\g(+e_2)$. This shows that all regular geodesics between two bases are homotopic.

Before proceeding, we show the following.\vs2

Claim 1: \emph{We may assume that $X$ is closed, i.e. a flat, in both $E_1$ and $E_2$.} \vs2

\emph{Proof:} Suppose not. Then $X$ is not closed in one of the matroids, say in $E_2$. Then the closure of $X$ in $E_2$ is $cl(X)=X+Z$. Let $E_3=E_2-Z$. Then $E_1\cap E_3=X$ but $|E_3|=r+k_3<r+k_2=|E_2|$. Since $(r-n,k_1+k_3)<(r-n,k_1+k_2)$, there is, by induction, a unique regular geodesic from $E_1$ to $E_3$ up to homotopy. Call this $\g_0$. The inclusion $E_3\subset E_2$ also gives a geodesic $\g_1$ up to homotopy. And it follows that any geodesic $\g$ from $E_1$ to $E_2$ is homotopic to $\g_0\g_1$. To see this delete $Z$ from each point in the path $\g$. This gives a geodesic $\g_2$ from $E_1$ to $E_3$ together with a homotopy to $\g\simeq \g_2$ which is fixed at the source $E_1$ and goes through a geodesic from $E_2$ to $E_3$ at the other end. This implies that $\g\simeq \g_2\g_1\simeq \g_0\g_1$ fixing the two endpoints proving the claim.\vs2

Now suppose that $n<r$ and $k_1+k_2>0$. By symmetry we may assume $k_1>0$. We also assume, by Claim 1, that $X$ is closed in $E_1$ and $E_2$.\vs2

Claim 2: \emph{If there exists a regular geodesic $\g$ from $E_1$ to $E_2$ which starts by adding an element of $E_2\backslash E_1$ to $E_1$, i.e. $\g(1)=E_1+e_2$, then any two regular geodesics from $E_1$ to $E_2$ are homotopic.}\vs2

\emph{Proof:} Since $X$ is closed in $E_2$, $X+e_2$ is independent for any $e_2\in E_2\backslash E_1$. So, by induction on $r-n$ there is a unique homotopy class of regular geodesics from $E_1+e_2$ to $E_2$. Let $\g(+e_2)$ denote the resulting geodesic from $E_1$ to $E_2$. Since $k_1>0$, $E_1$ contains a coindependent element $e_1$ which is not in $X$. Then $E_1-e_1\in IR(r,\FF_2)$ and, by induction on $k_1+k_2$, there is a unique homotopy class of regular geodesics from $E_1-e_1$ to $E_2$. Let $\g(-e_1)$ denote any of the resulting geodesics from $E_1$ to $E_2$. We claim that $\g(+e_2)\simeq\g(-e_1)$. This is because both geodesics can be chosen to pass through $E_1+e_2-e_1$ and the two 2-step path from $E_1$ to $E_1+e_2-e_1$ are homotopic since they are both contained in $C(E_1-e_1)$. Since this is true for any choice of $e_1$ and $e_2$, all regular geodesics are homotopic.\vs2

By Claim 2 we are reduced to the case when all regular geodesics from $E_1$ to $E_2$ start by deleting an element of $E_1$. This implies that $k_1\ge 2$ since, if $k_1=1$, $E_1+e_2$ will be regular for any $e_2\in E_2\backslash E_1$. Using the notation of the previous paragraph, we need to show that $\g(-e_1)\simeq\g(-e_2)$ for any $e_1,e_2\in E_1\backslash E_2$ which are coindependent in $E_1$. If $E_1-e_1-e_2$ has full rank then we will be done since we may then choose both path to pass through $E_1-e_1-e_2$. So, suppose not. Then $H=E_1-e_1-e_2$ is a hyperplane in $E_1$ which contains $X$. Since $k_1\ge2$, $H$ contains at least one coindependent element $e_0$ not in $X$. Then $E_1-e_0-e_i$ has full rank for $i=1,2$ which implies that $\g(-e_1)\simeq\g(-e_0)\simeq\g(-e_2)$ proving the last case of our lemma.
 \end{proof}
 
 \begin{lem}\label{second uniqueness lemma}
 If there exists a regular geodesic from $E_1$ to $E_2$ then any two regular geodesics are homotopic in $IR(r,\FF_2)$ (fixing the endpoints).
 \end{lem}
 
 
 \begin{proof}
 Let $X_0$ be a maximal independent subset of $X=E_1\cap E_2$ and let $Z=X-X_0$.
 Let $E_3=E_2-Z$. Then $E_1\cap E_3=X_0$ is independent. So, all regular geodesics from $E_1$ to $E_3$ are homotopic. If $\g_1,\g_2$ are regular geodesics from $E_1$ to $E_2$ and $\g_0$ is any geodesic from $E_2$ to $E_3$ then $\g_1\g_0$ and $\g_2\g_0$ are regular geodesics from $E_1$ to $E_3$. So, $\g_1\g_0\simeq \g_2\g_0$ which implies $\g_1\simeq\g_2$ as claimed.
 \end{proof}


\subsection{Uniqueness of paths up to homotopy}

 \begin{lem}\label{last lemma}
 If there exists a regular geodesic $\g_0$ from $E_1$ to $E_2$ then any two path from $E_1$ to $E_2$ are homotopic in $IR(r,\FF_2)$ (fixing the endpoints).
 \end{lem}

Taking $E_1=E_2$ this immediately implies \vs2

\noi{\bf Theorem \ref{thm:IR(F2) is simply connected}}
\emph{$IR(r,\FF_2)$ is simply connected.}

\begin{proof}[Proof of Lemma \ref{last lemma}]
By Lemma \ref{second uniqueness lemma}, any regular geodesics from $E_1$ to $E_2$ is homotopic to $\g_0$. So, it suffices to show that every path is homotopic to some geodesic. Let $\g$ be any path. Then the length of $\g$ can be written as
\[
    \el(\g)=d+\delta
\]
where $d=d(E_1,E_2)$ and $\delta\ge0$ will be called the \emph{excess} of $\g$. If $\delta=0$ then $\g$ is a geodesic and we are done. So, suppose $\delta>0$. We proceed by induction on the \emph{triple} $(\delta,|E_2|,d)$.

Claim: \emph{We may assume that $X=E_1\cap E_2$ is independent and closed in $E_2$}.

Otherwise, there is an element $x_0\in X$ which is coindependent in $E_2$. If we let $E_0=E_2-x_0$ and let $\g_1$ be the geodesic $E_2\to E_0$ then, first, $\g_0\g_1$ is a regular geodesic from $E_1$ to $E_0$ and, second, $\g\g_1$ is a path from $E_1$ to $E_0$ with triple $(\delta,|E_0|,d+1)<(\delta,|E_2|,d)$. Therefore, by induction on the triple, $\g\g_1\simeq\g_0\g_1$ which implies that $\g\simeq\g_0$ as required.


Now suppose that $\g$ is a path from $E_1$ to $E_2$ with excess-distance pair $(\delta,d)$ where $\delta>0$. Consider the first step $\g_1$ of $\g$ from $E_1$ to some $E_3$ and let $\g_2$ be the remainder. So $\g=\g_1\g_2$:
\[
    \g:E_1\xrarrow{\g_1}E_3\xrarrow{\g_2}E_2
\]
There are two possibilities: $d(E_3,E_2)=d\pm1$.

If $d(E_3,E_2)=d+1$ then $\g_1^{-1}\g_0$ is a regular geodesic from $E_3$ to $E_2$ and $\g_2$ has excess $\delta-2<\delta$. So, $\g_2\simeq \g_1^{-1}\g_0$ and $\g=\g_1\g_2\simeq \g_1\g_1^{-1}\g_0\simeq\g_0$ as required. So, we may assume that $d(E_3,E_2)=d-1$ (and $d>0$).

There are two possibilities. Either $E_3=E_1-e_1$ or $E_3=E_1+e_2$ and
\[
    E_3\cap E_2=\begin{cases} X & \text{if } E_3=E_1-e_1\\
    X+e_2 & \text{if } E_3=E_1+e_2
    \end{cases}
\]
In both cases $E_3\cap E_2$ will be independent by the Claim proved above. So, there exists a regular geodesic from $E_3$ to $E_2$ by Lemma \ref{existence lemma}. The path $\g_2$ has triple $(\delta,|E_2|,d-1)<(\delta,|E_2|,d)$. So, by induction on the triple, $\g_2$ is homotopic to a regular geodesic $\g_3$ from $E_3$ to $E_2$ which implies that $\g$ is homotopic to the regular geodesic $\g_1\g_3$ as claimed.

This concludes the proof of the lemma and the theorem.
\end{proof}
 
 Since $\pi_3 IR(3,\FF_2)=\ZZ^8$, the connectivity of $IR(r,\FF_2)$ is at most 2. But this leaves open the following question.
 
 \begin{quest}
 Is $IR(r,\FF_2)$ 2-connected for all $r$?
 \end{quest}
 



\end{document}